\documentclass[a4,11pt]{article}
\usepackage{amsmath,amsbsy,amsfonts,amssymb,graphics,epsfig,float,mathrsfs,amsthm,braket}

\usepackage{caption}
\usepackage{subcaption}
\usepackage{multirow}
\usepackage[retainorgcmds]{IEEEtrantools}

\usepackage{hyperref}
\hypersetup{colorlinks=true, linkcolor=blue}
\hypersetup{colorlinks=true,citecolor=red}

\renewcommand{\v}[1]{\mathbf{#1}}

\renewcommand{\c}[1]{\textsf{#1}}

\newcommand{\s}[1]{{\textsf{\textbf{#1}}}}
\newcommand{\pd}[2]{\frac{\partial #1}{\partial #2}}

\newcommand{\sd}[2]{\frac{\textrm{d} #1}{\textrm{d} #2}}
\newcommand{\sdd}[2]{\frac{\textrm{d}^2 #1}{\textrm{d} {#2}^2}}

\newcommand{\benE}[2]{\begin{IEEEeqnarray}{#1}\label{#2}}
\newcommand{\eenE}{\end{IEEEeqnarray}}
\newcommand{\be}[1]{\begin{equation}\label{#1}}
\newcommand{\ee}{\end{equation}}
\newcommand{\bes}[1]{\begin{equation}\label{#1} \begin{split}}
\newcommand{\ees}{\end{split} \end{equation}}
\newcommand{\ben}[1]{\begin{eqnarray}\label{#1}}
\newcommand{\een}{\end{eqnarray}}

\newcommand{\re}[1]{~(\ref{#1})}
\newcommand{\ret}[2]{~(\ref{#1},\ref{#2})}

\newcommand{\bma}{\begin{pmatrix}}
\newcommand{\ema}{\end{pmatrix}}

\newcommand{\vtu}{{\boldsymbol {\mathcal U}}}
\newcommand{\ctu}{{\mathcal U}}
\newcommand{\hctu}{{\widehat{\mathcal U}}}

\newtheorem{theo}{Theorem}
\newtheorem{prop}[theo]{Propostion}

\title{\huge\s{Design and stability of a family of deployable structures }}
\author{ \textsf{Thomas Lessinnes$^\dagger$ and Alain Goriely$^\ddagger$}\\
\\  {$\dagger$\it LCVMM, Ecole Polytechnique f\'ed\'erale de Lausanne, Switzerland} \\
{$\ddagger$ \it  Mathematical Institute, University of Oxford}\\
            {E-mail:\ {\tt
thomas.lessinnes\@@\,epfl.ch}}\\
            { \it Corresponding author: A.G.} }
\date{\today}

\begin{document}
\maketitle
\hrule\vskip 6pt
\noindent {\bf\small $\blacksquare$ Abstract\ } {\small
A large family of deployable filamentary structures can be built by connecting two elastic rods along their length. The resulting structure has interesting shapes that can be stabilized by  tuning the material properties of each rod. To model this structure and study its stability, we show that the equilibrium equations describing unloaded states can be derived from a variational principle. We then use a novel geometric method to study the stability of the resulting equilibria. As an example we apply the theory to establish the stability of all possible equilibria of the Bristol ladder. 
 }\vskip 6pt
\hrule


\section{Introduction}

\begin{figure}[h]
\centering
\includegraphics[width=15cm]{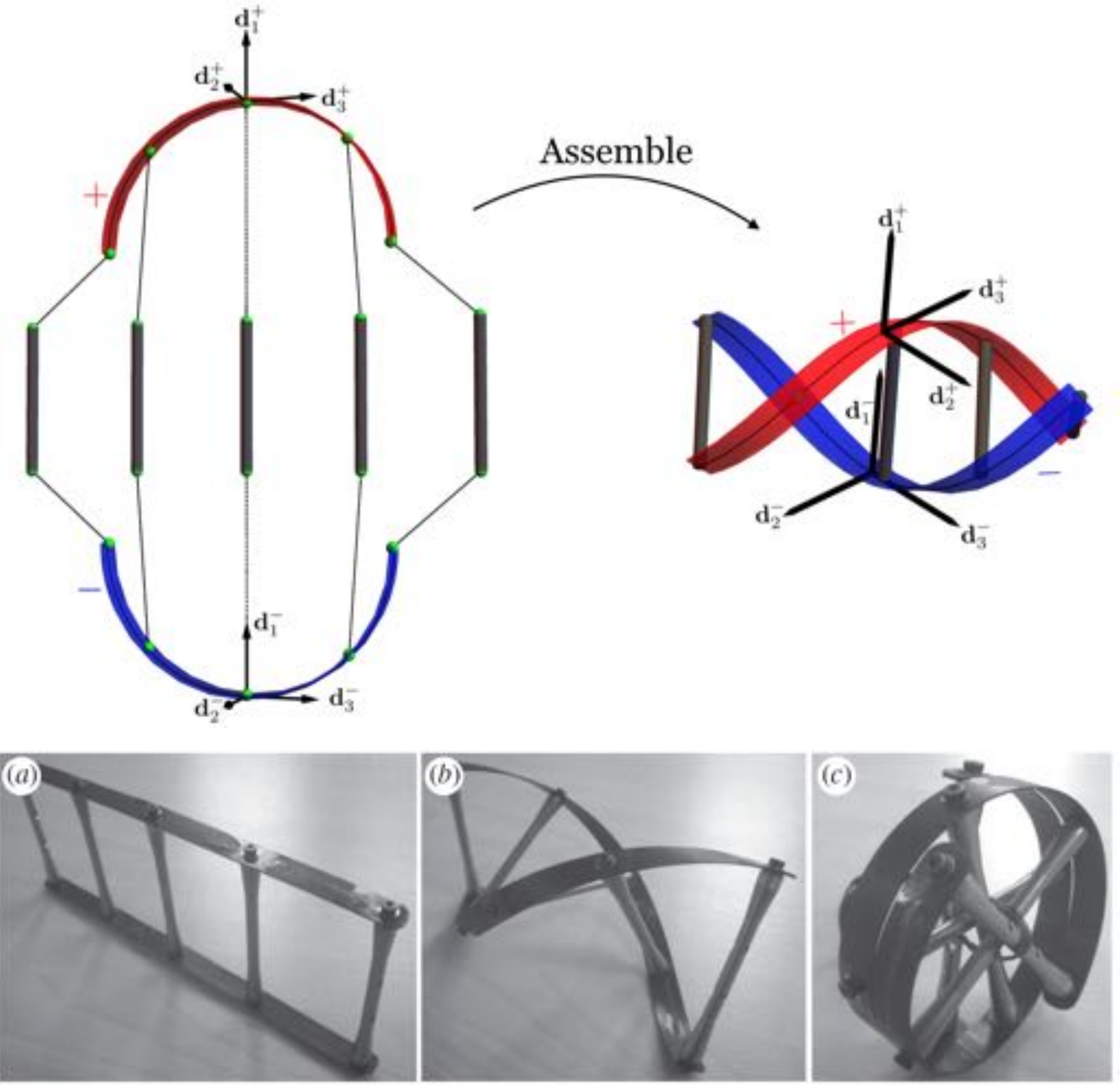}
\caption{Top: A sketch of the Bristol ladder consisting of two flanges connected by rigid spokes. Bottom: Actual ladder as developed in Paul Weaver's laboratory (From~\cite{laweda12} reproduced without authorization)} \label{fig-ladder}
\end{figure}

The \textit{Bristol ladder}~\cite{laweda12} is an ingenious structure composed of  two pre-curved flanges with rectangular cross-sections connected by rigid spokes (see Figure~\ref{fig-ladder}). This structure can be almost completely folded into itself by rotating the end, it can also be deployed to act as a tether or to connect other larger structures. Remarkably, it has multiple stable helical and straight configurations that can be activated under different loads. 

We use this ladder as a general motivation to consider the theoretical stability and design of a more general class of \textit{deployable ladders} composed of two uniform unshearable Kirchhoff rods constrained by rigid spokes maintaining constant distance along the length of the structure. Following the original work, we assume that the flanges are allowed to rotate freely about the axes of the spokes. In our generalization, we allow each structure to have arbitrary intrinsic curvature.
The equations governing the equilibria are obtained as the states for which the first variation of the elastic energy functional vanishes for all admissible perturbations. We show that these extra degrees of freedom introduced in our model can be tuned to stabilize the structure and create long filaments with new interesting shapes and properties. Of particular importance for the context of this paper,  we show that specific designs can be achieved through a detailed rigorous stability analysis of the different equilibrium states.

The notion of stability used here is based on the second variation of the elastic energy functional being strongly positive-definite at stable equilibria. In a recent paper, we showed that the study of the second variation for one-dimensional problems deriving from an energy potential  can be reduced to the geometric study of orbits in the phase plane associated with such a potential \cite{lego15}. Here, we make explicit use of these results to design multi-stable ladders by controlling particular orbits of interest and assessing directly their stability properties through their geometry: without having to consider or even compute the second variation of the energy functional. In particular, we consider the possibility of perversion consisting of a helical structure  with one handedness connected to a  helical structure with the reverse handedness \cite{gota98b,mcgo02}.
\section{The structure} 

\subsection{Shape and geometry}\label{sec-geom}

The deployable ladders considered here are a special case of a class of  birods that we studied previously~\cite{lemo15}. The description of these ladders is based on the general theory of Kirchhoff rods~(see \cite{molego12} for an introduction and terminology).  Each flange is labelled by a sign $-$ (minus) or $+$ (plus) and all associated quantities take the corresponding superscript. A flange is modeled as a spatial curve $\v r^\pm$, the centerline of the flange, and is equipped with a family of local material frames of directors $\{ \v d_1^\pm, \v d_2^\pm, \v d_3^\pm\}$, where $\v d_3^\pm$ is tangent to the centerline of the $\pm$ flange, $\v d_1^\pm$ and $\v d_2^\pm$ span its cross-sections and $\v d_1^\pm$ is normal to its tangent plane (cf. Figure~\ref{fig-ladder}). These quantities are functions of the reference arc-length $R$  of the minus flange expressed in units of the length $a$ of the spokes. That is, $R$  is the arc length along the centerline when the flange assumes zero stress state (the configuration shown on the top left of Figure~\ref{fig-ladder}). The shape of a flange is then fully described by a scalar function $\alpha^\pm$, the stretch of the flange, and a vector function $\v u^\pm$, called the Darboux vector of the family of frames $\{ \v d_1^\pm,\, \v d_2^\pm,\, \v d_3^\pm\}$ such that 
\begin{align}
\label{shape1rod}
\left( \v r^\pm \right)'&= \alpha^\pm ~ \v d_3^\pm,   
&
\left( \v d_i^\pm\right)' &= \v u^\pm \times  \v d_i^\pm,   
\end{align}
where $(~)'$ denotes derivatives w.r.t. $R$. Note that $\alpha$, $\v u$ and Equations~(\ref{shape1rod}) encode the shape of a flange up to a global translation and rotation that can be specified by providing initial conditions for\re{shape1rod}. In any given configuration of a flange, the components $\c u_i$ of $\v u$ in the local basis provide the rate of bending of the flange per unit of $R$ about the directions $\v d_i$.

In the reference state of each flange, that is the stress free state showed on the top left panel of Figure~\ref{fig-ladder}, the stretch is $\widehat{\alpha}^\pm=1$. For the Bristol ladder, the flanges are prepared so that they are naturally curved about the $\v d_2^\pm$ direction: their reference Darboux vectors are $\widehat{\v u}^\pm= \mp\, \widehat {\c u}~\v d_2^\pm$, where $\widehat {\c u}$ is a positive constant. For a more general deployable ladder, we also consider the possibility of pre-curving the flanges about the $\v d_1^\pm$ direction. Therefore, we assume that $\widehat{\v u}^\pm=  \widehat {\c u}_1^\pm \v d_1^\pm + \, \widehat {\c u}^\pm~\v d_2^\pm$ where all four quantities $\widehat {\c u}_1^\pm$ and $\widehat {\c u}^\pm$ are independent of $R$ (see Appendix~\ref{app-energy} for the general case). These four quantities only appear in the problem in the following three linear combinations 
\begin{align}
\label{defU}
\hctu_0&= \widehat{\c u}_1^+ - \widehat{\c u}_1^-,
&
\hctu_1&= \widehat{\c u}_1^+ + \widehat{\c u}_1^-,
& 
\widehat u&= \frac{\widehat{\c u}_2^- - \widehat{\c u}_2^+}2.
\end{align}
For the Bristol ladder, we have $\hctu_0= \hctu_1= 0$ and $\widehat u=1$.

When the structure is assembled, the rigid spokes insure that $\v d_1^+(R) = \v d_1^-(R)=\v d_1(R)$ for all $R$. A configuration of the structure is then fully specified by four scalar functions of $R$: $\alpha^+$ and $\alpha^-$ the stretches of each flange, $\theta$ the angle between $\v d_3^-$ and $\v d_3^+$, and $\ctu $ the rate of bending of the structure about the direction $\v d_1$~\cite{lemo15} (see also Appendix~\ref{app-birod}). We  study here the inextensible case  ($\alpha^+ = \alpha^-=1$) relevant for structures such as the Bristol ladder.

\subsection{Internal energy}

The starting point for computing the equilibrium states of such a structure is to compute its elastic energy. Based on the elastic energy of each flange and their geometric constraints, we show in Appendix~\ref{app-energy} that the internal energy of the ladder in a state $(\theta(R), \ctu(R))$ is given by
\ben{LadEn}
\mathscr E= \frac{ (EI)}{4} \int_0^{L}   
\frac{1}{2}\left (\theta'-\hctu_0\right )^2 + 2 \left ( \, \ctu-\frac{\hctu_1}{2 }\right )^2- \bigg ( 4 b \cos \theta ~\left( 1 - \widehat{ u} \right) - (b-\Gamma) \,  \cos 2 \theta \bigg)~ \textrm{d} R,
\een
where all lengths are scaled in units of $a$ the length of the rigid spokes, $(EI)$ is the bending stiffness of the flanges about the $\v d_1$ direction, $b$ is the ratio of bending stiffnesses about the $\v d_2$ and $\v d_1$ direction, $\Gamma$ is the ratio between the torsional stiffness of the flanges, and $L$ is the length of the flanges. As shown explicitly in Appendix~\ref{app-originalmodel}, the original  model~\cite{laweda12} for the Bristol ladder is recovered by imposing the constraint $\theta'=\ctu=0$ together with the choice of parameters $\widehat u =1$ and $\hctu_0 =\hctu_1=0$. These constraints are obtained in the limit $b\to 0$ (which implies $\Gamma\to 0$ see Appendix) so that the first two terms in\re{LadEn} must be minimized independently.
As for the choice of parameters, we first consider the case $\widehat u= 1$ as it allows to best display our method. The general case for arbitrary $\widehat u$ is briefly discussed in Section~\ref{sec-effectofu}.

Stable equilibria of the structure  are local minima of the energy functional\re{LadEn}. If the pair $(\theta^\star, \, \ctu^\star)$ is a local minimum of the energy, then $\ctu^\star = \hctu_1/2$ since any other function $\ctu$ would increase the second term in\re{LadEn} while leaving all other terms unchanged. Therefore, on stable equilibria, we can  further simplify our model by setting $\ctu= \hctu_1/2$ so that the second term in\re{LadEn} disappears. 
By scaling the energy w.r.t. $(EI)/4$ and defining the positive constant $\varepsilon = b-\Gamma>0$, the energy of the structure takes the non-dimensional form
\be{LadEnnd}
\mathscr E_{\textrm{n.d.}} =  \int_0^{L}   {\mathscr L(\theta,\theta') ~\textrm{d} R=\int_0^{L} \frac{1}{2}\left (\theta'-\hctu_0\right )^2 +  \varepsilon \,  \cos 2 \theta }~ \textrm{d} R.
\ee
Equilibria of the ladder are solutions of the Euler-Lagrange equations  associated with\re{LadEnnd} with free-boundary conditions: 
\be{ELlad0}
-\left(\frac{\partial \mathscr L}{\partial\theta'}\right )' + \frac{\partial \mathscr L}{\partial\theta}=0, \qquad \textrm {and} \qquad 
\left. \pd {\mathscr L}{\theta'} \right |_{R=0} =\left. \pd {\mathscr L}{\theta'} \right |_{R=L} = 0.
\ee

\section{Qualitative analysis}

An equilibrium is defined as \textit{stable} if it locally minimizes the energy of the system. A sufficient condition for stability is that the second variation of the energy functional be strongly positive for all admissible perturbations (see e.g.~\cite{gefo00,Manning:2009cq}). Conversely, if there exists an admissible perturbation such that the second variation of the potential energy is strictly negative, then the equilibrium is \textit{unstable}. 
 
Our system supports \textit{constant equilibria} (i.e. solutions with constant $\theta$) as well as \textit{non-constant equilibria}. Stability for equilibria can  be established by studying the positive-definiteness of the second variation of the functional\re{LadEnnd}. For constant equilibria, this analysis is straightforward. However, for non-constant equilibria, establishing positive-definiteness is typically intractable analytically. Recently, we established that in the case of one-dimensional functionals such as\re{LadEnnd}, a geometric analysis of the trajectory in the phase plane is sufficient to establish stability and instability~\cite{lego15}. In this  framework, stability can be established without explicit knowledge of the analytical form of the solution.

 We start by summarizing the key results. Then we perform a qualitative analysis of the possible equilibria and their stability. 
 We complete our analysis with a  quantitative analysis of all possible equilibria in Section~\ref{sec-quantit}.

\subsection{Stability of one-dimensional systems: a summary\label{sec-sta1D}}
 A functional  of the form, 
\ben{myfunc}
\mathscr E[\theta]=  \int_a^{b} \mathscr L\big( \theta(s),\theta'(s) \big)   \textrm{d} s,
\een
with
\be{FormL}
\mathscr L[ \theta,\theta']= \frac{(\theta'-A)^2} 2 - V(\theta),
\ee
where $V$ is a $C^2$ function such that its first and second derivative never vanish simultaneously and $A$ is a real constant,  has a Euler-Lagrange equation of the form 
\begin{equation}\label{thetapppot1}
\theta''(s) + \left .\sd{V}{\theta}\right|_{\theta(s)} = 0.
\end{equation}

In the phase plane  of Equation\re{thetapppot1}, we define the sets $\Gamma_M=\{(\theta,\theta')\in\mathbb R^2\big|~ \theta'\in \mathbb R,\, V'(\theta)=0, \,V''(\theta)<0\}$ and  $\Gamma_m=\{(\theta,\theta')\in\mathbb R^2\big |~ \theta'\in \mathbb R, \,V'(\theta)=0,\, V''(\theta)>0\}$ which are the vertical lines  at maximal and minimal points of $V$. Given a solution $\theta^\star$ of\re{thetapppot1} with domain $[s_1,s_2]$, we define the associated trajectory $\eta$ in the phase plane: 
\be{defeta}
\eta:[s_1,s_2]\to \mathbb R^2:s\to (\theta^\star(s),{\theta^\star}'(s)).
\ee
 Finally, for any such trajectory $\eta$, we define the \textit{stability index}
\be{defJ}
J[\eta]= \# \left \{ s\in[s_1,s_2]: \eta(s)\in \Gamma_m\right\}-\# \left \{ s\in[s_1,s_2]: \eta(s)\in \Gamma_M \right\},
\ee
as the number of times $\eta$ crosses minimal points of $V$ minus the number of times it crosses maximal points of $V$. 

Assuming that we are interested in states $\theta$ which minimise $\mathscr E$ while respecting natural boundary condtions $\left . \pd{\mathscr L}{\theta'}\right|_{\theta(a)}=\left . \pd{\mathscr L}{\theta'}\right|_{\theta(b)}=0$,   the stability of $\theta$ is given by the following result~\cite{lego15}:
\begin{theo} \label{th-indJ}
Let $\theta:[a,b]\to \mathbb R$ be a stationary function of the functional $\mathscr E$  such that $\exists s\in[a,b]:\theta'(s)\neq 0$ and $\nexists s\in[a,b]: \theta'(s)=0~ \& \left . \sd{V}{\theta}\right|_{\theta(s)}=0$. Let $\gamma:s\in[a,b]\to (\theta(s),\theta'(s))\in\mathbb R^2$ be the trajectory associated to $\theta$ in the phase plane.

If $J[\gamma]>0$, then $\theta$ is not a local minimum of $\mathscr E$. If $J[\gamma]<0$, then $\theta$ is a local minimum of $\mathscr E$.
\end{theo}

 If $A=0$ some solutions of\re{thetapppot1} may be constant functions. Such functions do not satisfy the hypothesis of Theorem~\ref{th-indJ}. In that case, their stability can be readily established (see e.g. \cite{lego15}):
\begin{prop} \label{prop-sta}
Let $\theta:[a,b]\to \mathbb R$ be a stationary function of the functional $\mathscr E$  such that $\forall s\in [a,b]: \theta'(s)=0$,  then  $\theta$ is a local minimum of $\mathscr E$ if and only if $\left. \sdd{V}{\theta}\right |_\theta<0.$
\end{prop}

\subsection{Application to the deployable ladder} 

The energy functional\re{LadEnnd} of the ladder has the form\ret{myfunc}{FormL} with 
\be{defV}
V(\theta) = - \varepsilon \, \cos 2\theta.
\ee
The equilibria of the ladder are defined as the solutions of Equation~(\ref{ELlad0}) given explicitly by
\be{ELlad}
\theta'' +2\, \varepsilon \sin 2 \theta =0, \quad \textrm{with boundary conditions\quad }\theta'(0)=\theta'(L)=\hctu_0.
\ee

We start our analysis with the special case $\hctu_0=0$. Then, the boundary value problem\re{ELlad} admits three constant solutions $\theta_{\textrm{st}}(R)=0$, $\theta_{\textrm{r}}(R)=-\pi/2$, and $\theta_{\textrm{l}}(R)=\pi/2$ which respectively correspond to straight, right-handed  and left-handed helical equilibria. The stability of these equilibria is obtained from Proposition~\ref{prop-sta}. We have 
\begin{equation}
\left . \sdd{V}{\theta}\right |_{\theta_\text{r}} = \left . \sdd{V}{\theta}\right |_{\theta_\text{l}} = - 4 \varepsilon < 0, \quad \textrm{and \quad }\left . \sdd{V}{\theta}\right |_{\theta_\text{st}} = 4 \varepsilon >0,
\end{equation}
 so that the helical solutions $\theta_\text{r}$ and $\theta_\text{l}$ are stable while the straight solution $\theta_\text{st}$ is unstable.  

\begin{figure}[h]
\centering
\includegraphics[width=12cm]{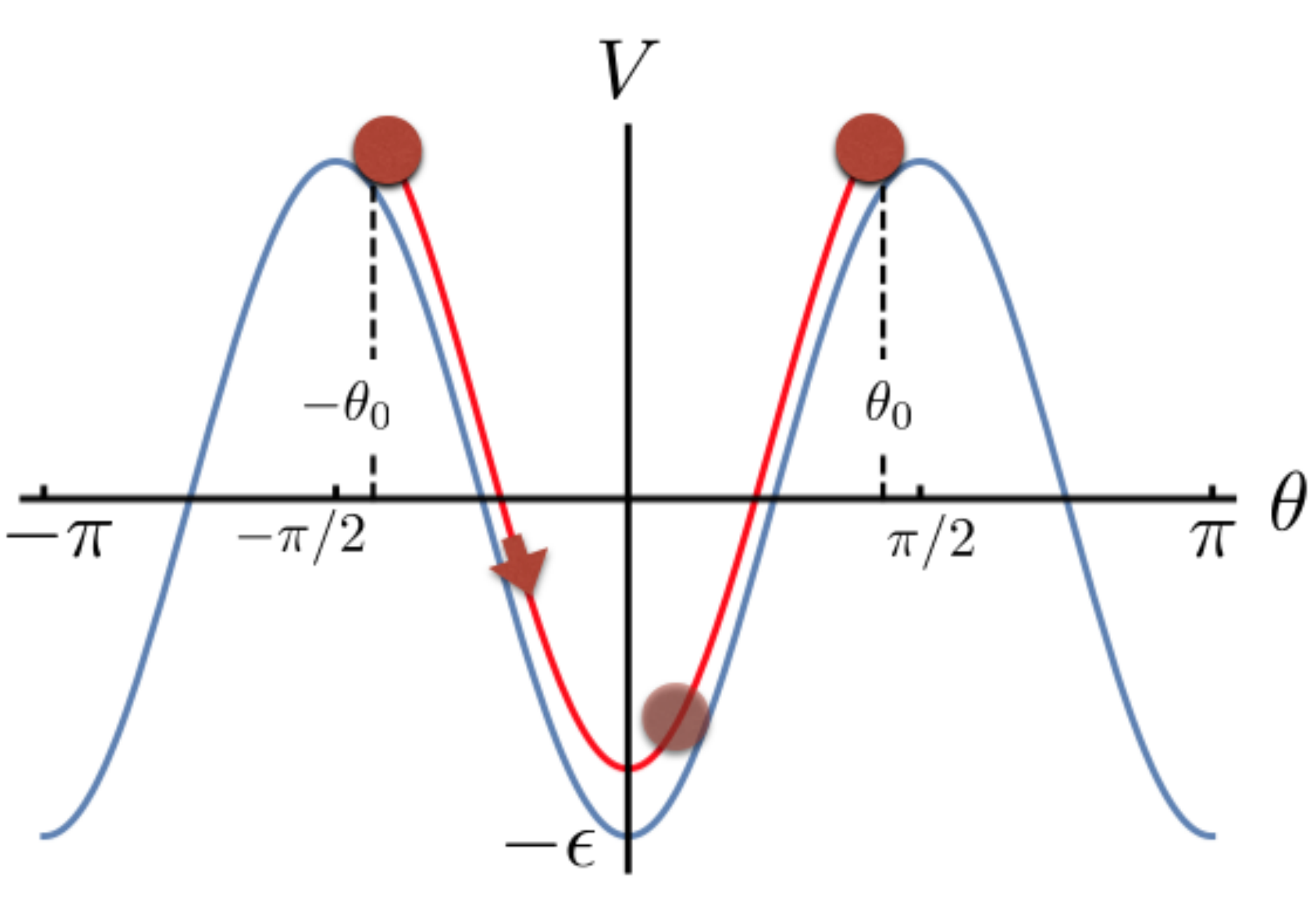}
\caption{Equation\re{ELlad} is equivalent to that of a point mass particle submitted to the potential $V=-\varepsilon \cos 2 \theta$. The boundary value problem therefore amount to finding all possible trajectories of the particle such that it starts and ends with horizontal velocity $\hctu_0$ and such that the flight time is $L$.} \label{fig-slidonV}
\end{figure}

Depending on the length $L$ of the structure, non-constant solutions of\re{ELlad} may also exist. Indeed, by analogy to a point-mass mechanical system, the differential equation\re{ELlad} also describes the motion of a point mass  in a potential energy $V$ (see Figure~\ref{fig-slidonV}). Hence, our system admits the \textit{pseudo-energy} $E$: 
\be{ConsE}
E=\frac{\theta'^2}{2}+ V(\theta). 
\ee

Given a number $\theta_0\in (0,\pi/2)$, there exists a ladder of length $L= L_p(\theta_0)$ which admits an equilibrium corresponding to the solution of\re{ELlad} shown in Figure~\ref{fig-slidonV}: it starts with $\theta(0)=-\theta_0$, swings through the well of potential at $\theta=0$ and ends at $\theta(L)=\theta_0$. The function $L_p$ is defined by
\be{Lth0}
L_p(\theta_0)= \int_0^L \textrm{d} R= \int_{-\theta_0}^{\theta_0} \frac{\textrm{d}\theta}{\theta'} = \int _{-\theta_0}^{\theta_0} \frac{\textrm{d}\theta}   {\sqrt{2 (V_0 - V(\theta))}},
\ee
where $V_0=V(\theta_0)$ and the last equality comes after substituting $\theta'$ according to\re{ConsE}.

In the limit $\theta_0\to \pi/2$, the integral\re{Lth0} diverges so that $\lim_{\theta_0\to \pi/2} L_p(\theta_0)=+\infty$. In the limit $\theta_0\to 0$, the potential can be approximated by its osculating parabola: $V(\theta)= \varepsilon \left (-1 + 2 \theta^2 +O(\theta_0^4)\right )$. After substituting $V$ for $V_{\textrm{osc}}= - \varepsilon (1-  2 \theta^2)$ in\re{Lth0}, we find $\lim_{\theta_0\to 0} L_p(\theta_0)= \frac{\pi}{2 \sqrt{\varepsilon}}=:L_0$.  Since $L_p(\theta_0)$ is a continuous function, the solution described in Figure~\ref{fig-slidonV} exists at least for all ladders of length $L\in(L_0,+\infty)$. Since these solutions invert their handedness we refer to them as \textit{perversion}. Following the same type of argument, it is easy to check that the condition $L\geq L_0$ is necessary and sufficient as there  is no such solution with $L<L_0$.

Finally, we note that there also exist equilibria with multiple passages through $\theta=0$ (multiple perversions). For a ladder of length $L$ such that there exists a natural number $n \in \mathbb N_0$ with $ n \leq L/L_0 < n+1$, there exists exactly one equilibrium with $k$ perversions for each odd number $k\leq n$ and exactly two equilibria with $j$ perversions for each even number $j\leq n$. 

All perversions and multiple perversion equilibria are unstable. Indeed, any solution $\theta$ of\re{ELlad} with $k$ oscillations  (where $0< k\in\mathbb N_0$) in its potential well crosses $k$ times the minimum of $V$ at $\theta=0$ and never crosses any of the maxima of $V$. The associated phase plane trajectory $\gamma:[0,L]\to \mathbb R^2:R\to (\theta(R),\theta'(R))$, is such that $J[\gamma]=k>0$  and Theorem~\ref{th-indJ} implies that $\theta$ is not a local minimum of the functional $\mathscr E_{\textrm{n.d.}}$. 

\subsection{Designing a stable perversion\label{sec-stabper}} 

From an engineering point of view, the presence of perversions could either be catastrophic occurrences -- e.g. hindering the proper deployment of a satellite arm -- or can be useful. Indeed, perversions have been shown to behave as twistless spring, that is their response in extension is much closer to an ideal spring than a single helical rod that would lock-up at a certain extension \cite{mcgo02}. It is therefore of interest to control the stability of a perversion, be it to enforce it or to avoid it.

\begin{figure}[h]
\centering
\includegraphics[width=12cm]{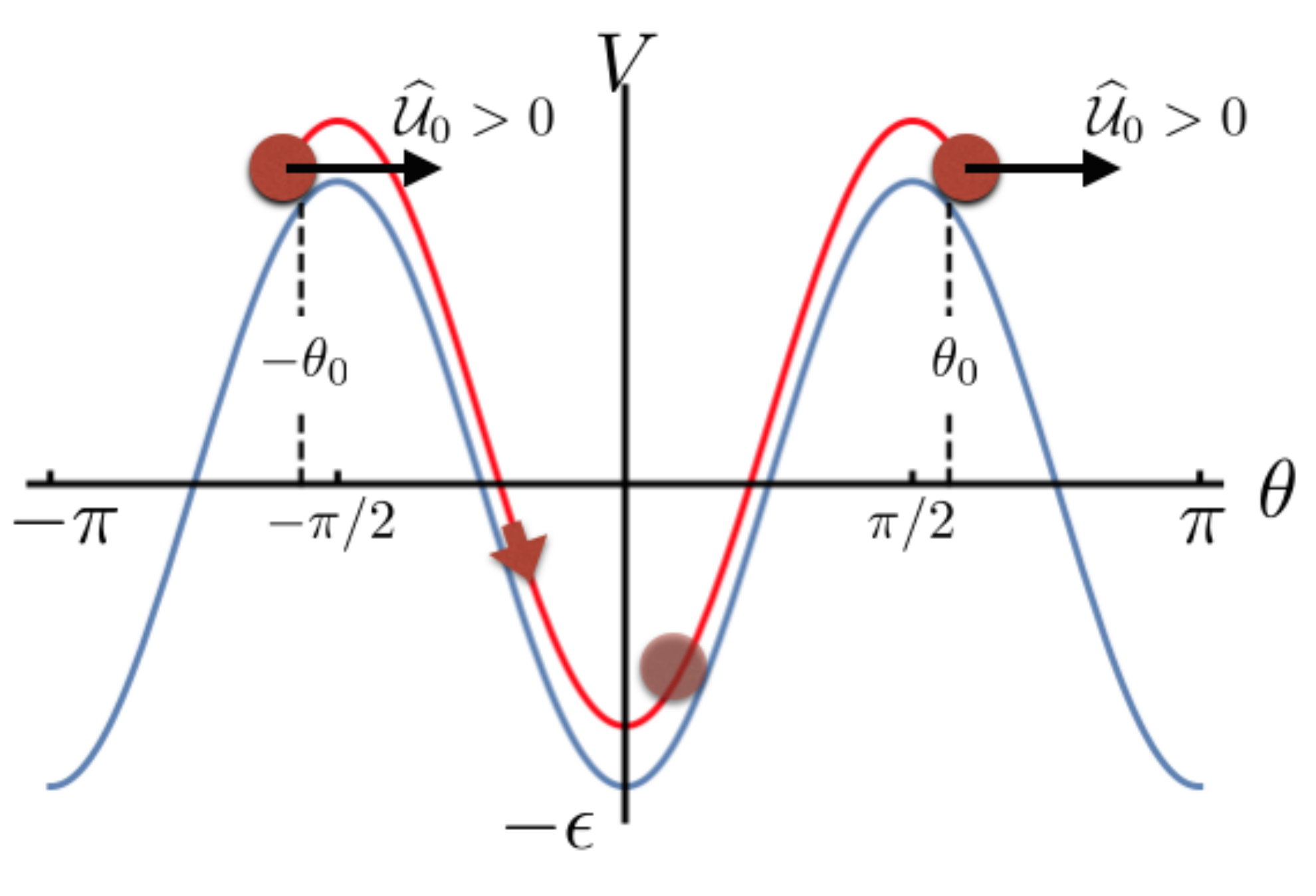}
\caption{A solution of Equation\re{ELlad} (with $\hctu_0>0$) that crosses both maxima of $V$. The index of the associated phase plane trajectory is $J[\gamma]=1-2=-1$ and this solution is therefore a  stable equilibrium of the ladder. } \label{fig-stabpervert}
\end{figure}

Perversions are unstable if $J[\gamma]>0$  (Theorem~\ref{th-indJ}). Since by definition, a perversion crosses the minimum of $V$ at $\theta=0$, a stable perversion  must also cross at least one of the maxima at $\theta=\pm \pi/2$. The analogy is helpful once again: to cross a maximum of $V$, the point mass needs initial kinetic energy: $\hctu_0> 0$ as shown in Figure~\ref{fig-stabpervert}. Recalling, from \re{defU} that $\hctu_0= \widehat{\c u}_1^+ -\widehat{\c u}_1^-$,  $\hctu_0$ can be made positive by warping the flanges so that  $\widehat{\c u}_1^\pm\neq 0$. If we further impose that the axis of the ladder remains straight, we require $\ctu= \hctu_1/2=0$ or more directly $\widehat{\c u}_1^-=- \widehat{\c u}_1^+$. 

\begin{figure}[h]
\centering
\includegraphics[width=15cm]{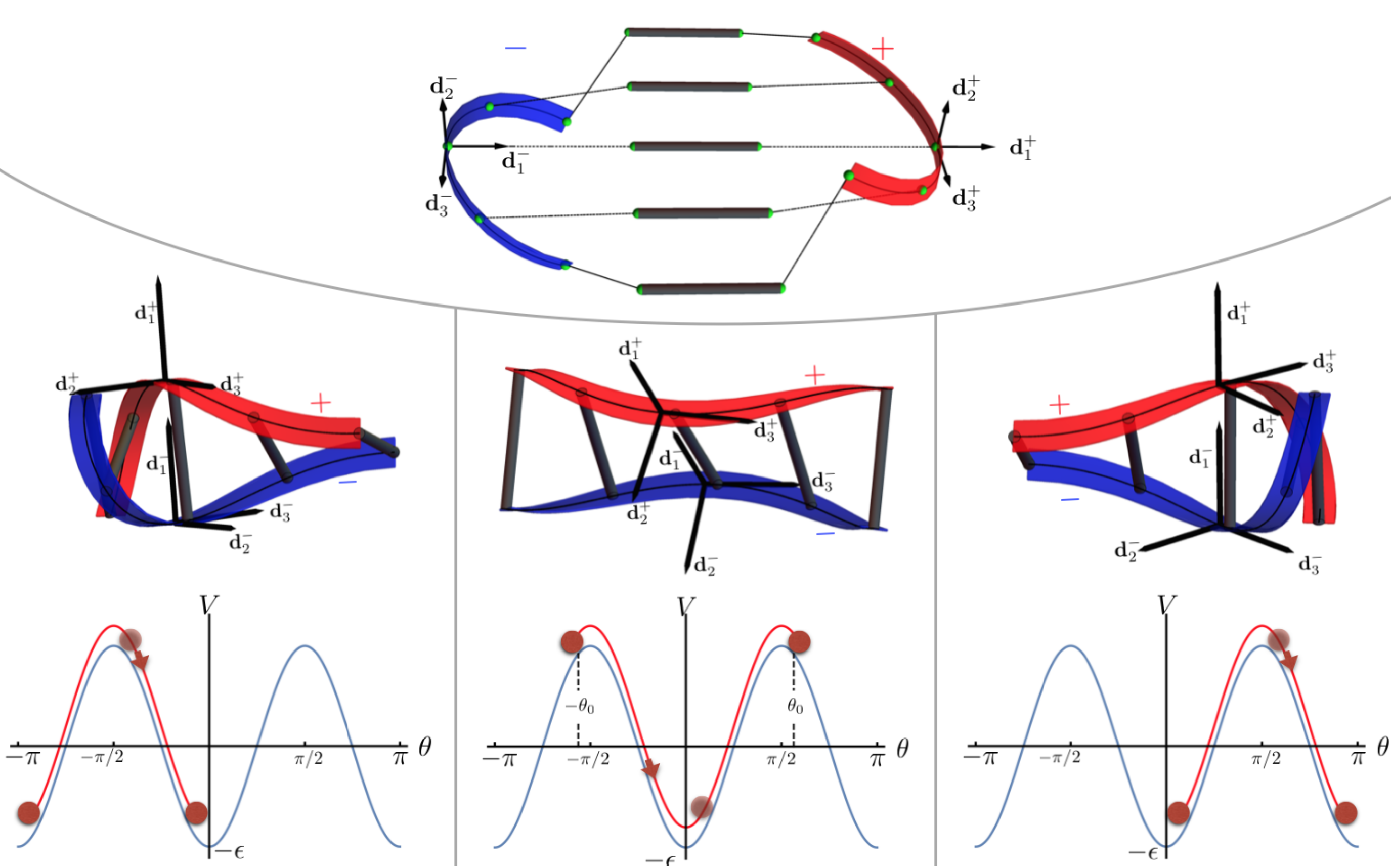}
\caption{A tri-stable ladder. The top panel show the flanges in their reference states with $\widehat{\v u}^\pm =  \pm \left ( \frac{\widehat {\c u}}{2} \v d_1^\pm +\widehat{\c u }\, \v d_2^\pm \right )$ -- where for this figure we took $\widehat{\c u} = 1/a$. Their reference centrelines lie on circles of radius $\frac 2{\sqrt{5}} \widehat{\c u}^{-1}$.  The three bottom panels show (all) the stable equilibria of the resulting ladder.   \label{fig-ladder1}}
\end{figure}

 If we assume that  $\hctu_0\gg \sqrt{2 \varepsilon}$,
 the total energy of a point mass in the potential well is dominated by its kinetic part and its speed is approximatively constant: $\theta'=\hctu_0$. It travels from $-\pi/2$ to $\pi/2$ during a period of time $L\simeq \pi/\hctu_0$. 
  For a ladder of length $L=\pi$, a stable perversion requires $\hctu_0\gtrsim 1$, that is, the flanges need to be pre-curved at least half as much about $\v d_1^\pm$ as they were about $\v d_2^\pm$:  $\widehat{\c u}_1^\pm\simeq \pm   a^{-1}/2 =  \widehat{\c u}/2$. 
A structure designed with such specifications would be tri-stable. A specific example of such a ladder is shown in  Figure~\ref{fig-ladder1} with its three stable configurations. Note that if $\hctu_0$ is too large, the solutions of single handedness will cross the minima of $V$ at $0$ and $\pm \pi$ which will destabilise them since $J[\gamma]$ would switch from -1 to $1$. Accordingly there exists a range of values of $\hctu_0$ for which the ladder is tri-stable. If $\hctu_0$ is too small, only the right and left solutions are stable. If $\hctu_0$ is too large only the perversion is stable. In the next section, we compute this range explicitly.
 
\section{Quantitative analysis \label{sec-quantit}} 

The qualitative analysis provided insight into the existence of equilibria with given stability properties. This analysis can be complemented by a quantitative analysis of these solutions by solving explicitly the  boundary value problem\re{ELlad}. 

\subsection{Exact solution for perversions\label{sec-quantper}}
We define a perversion of length $L$  as any solution of the differential equation\re{ELlad} such that both $\theta\big(L-R\big) = - \theta\left(R \right )$ and $\theta'(R)\neq0~\big(\forall R\in(0,L)\big)$. In particular, we do not enforce the boundary conditions $\theta'(0)=\theta'(L)=\hctu_0$.

The perversions can be characterized by their length $L$, their amplitude $\theta_0=\theta(L)$ and their non-dimensional pseudo-energy through $\nu=(E/\varepsilon-1)/2$ (with $\nu\in[-1,+\infty)$) where $E$ was defined in\re{ConsE}.

Inverting Equation~(\ref{ConsE}), leads to 
\be{thetap}
\theta' = \textrm{Sign}(\theta_0) \sqrt{ 2 \big(E- V(\theta)\big)} =  \textrm{Sign} (\theta_0)~ 2 \sqrt{\varepsilon({1+\nu})\left(1-\frac{1}{1+\nu} \sin^2 \theta\right)}.
\ee
 And we conclude that a perversion is determined by its length $L$ and the non-dimensional parameter $\nu$. 

The first order differential equation\re{thetap} with the initial value $\theta(L/2) =0$   admits a solution given implicitly by:
\be{solthetap}
 \sqrt \varepsilon L~\sqrt{1+\nu} \left (\frac {2 R}{L}-1\right )  = \textrm{Sign}(\theta_0)~F\left (\theta(R)\left | \frac{1}{1+\nu}\right. \right ) ,
\ee
where $F(\varphi | m)$
denotes the incomplete elliptic integral of the first kind. 

Equation\re{solthetap} can be used to quantify the existence of stable perversions. First note that a perversion $\theta$ is a monotonic function of $R$: when $\theta_0>0$ the left-hand side of\re{solthetap} increases by $F(\frac \pi 2 |\frac{1}{1+\nu})$ in every $\pi/2$ increment of  $\theta$. Furthermore, this solution also provides a relation between  the three parameters:
\be{th0pre}
\sqrt{\varepsilon} L \sqrt{1+\nu} = \textrm{Sign} (\theta_0) F\left( \theta_0 \left | \frac{1}{1+\nu} \right . \right ).
\ee
This relation provides an explicit expression of $|\theta_0|$ as a function of $\nu$: 
\be{th0} 
|\theta_0|= \arccos \left [ \textrm{cn} \left ( (\sqrt \varepsilon\, L)~\sqrt{1+\nu} \left |\frac{1}{1+\nu}\right.   \right ) \right ],
\ee
where $ \textrm{cn}(u|m)$ is the Jacobi elliptic cosine with parameter $m$.

Since the ladder is a physical object, it cannot interpenetrate. A necessary condition for the absence of hard contact between the flanges  is $|\theta_0|<\pi$ or equivalently 
\be{nohardcontact}
\sqrt{\varepsilon} L\, \sqrt{1+\nu}<2 K\left(\frac{1}{1+\nu} \right),
\ee
where $K(m) = F(\frac \pi 2 |m)$ is the complete elliptic integral of the first kind.

So far, we have established global properties of the perversion as a solution of the differential equation\re{ELlad} independently of the boundary conditions.  For a perversion to be an equilibrium of the ladder, it must also respect the boundary values: $\theta'(0)= \theta'(L)= \hctu_0$. Substituting this condition in\re{ConsE} leads to
\be{ateq0}
E= \varepsilon (1+ 2 \nu) = \frac{\hctu_0^2}{2} + V(\theta_0) = \frac{\hctu_0^2}{2} -\varepsilon \cos 2 \theta_0. 
\ee
Finally, substituting\re{th0} in\re{ateq0} yields a necessary and sufficient condition on $\nu$ for a perversion to be an equilibrium of the system:  
\be{Eqper}
\nu + \textrm{cn}^2\left (\sqrt{\varepsilon}\, L ~\sqrt{1+\nu} \left | \frac{1}{1+\nu} \right. \right)=\frac {\hctu_0^2}{4 \varepsilon}. 
\ee
In particular, (only) two constitutive parameters are relevant: $\ell=\sqrt{\varepsilon}\, L$ and $y=\hctu_0^2/(4 \varepsilon)$. 

We note that if $|\theta_0|<\pi/2$, a perversion crosses only one extremum of $V$: the minimum at $\theta=0$. In this case, we have $J[\gamma]=+1$ and we conclude that such a perversion is unstable (cf. Theorem~\ref{th-indJ}). If $|\theta_0|>\pi/2$, the perversions   cross both maxima at $\theta= \pm \pi/2$ and we have $J[\gamma]=-1$. These perversions are stable as local minima of the energy functional. Using\re{solthetap}, the condition $|\theta_0|>\pi/2$ becomes 
\be{stasol}
\sqrt{\varepsilon} L\, \sqrt{1+\nu}> K\left(\frac{1}{1+\nu} \right).
\ee
Conditions\ret{nohardcontact}{stasol} can be written as
\be{ineqper}
\frac{\sqrt{\varepsilon} L}{2} <k(\nu)<\sqrt{\varepsilon} L,\quad\text{where}\quad k(\nu)=\frac{1}{\sqrt{1+\nu}}K\left(\frac{1}{1+\nu} \right).
\ee
The function $k(\nu)$ is monotonically decreasing (see Figure~\ref{fig-k}). Accordingly, the inequalities\re{ineqper} determine an interval of admissible values of $\nu$ that guarantees the existence of stable perversions: the upper bound corresponds to the requirement that there should be no hard contacts between the flanges and the lower bound guarantees their stability. 

We conclude that a ladder of length $\ell$ (expressed in units of $a/\sqrt{\varepsilon}$) admits stable perversions for each value of $\nu$  that solves Equation\re{Eqper} in the interval  defined by\re{ineqper}.

Since $\theta(R)\leq \theta_0<\pi$ $\forall R\in[0,L]$, the solution can be written explicitly as
\be{}
\theta(R)=\textrm{Sign}\left(R-\frac{L}{2}\right) \textrm{arccos} \left [ \textrm{cn} \left(
\sqrt{\varepsilon}L\, \sqrt{1+\nu} \left(2R/L-1\right) \left| \frac{1}{1+\nu} \right.\right)\right].
\ee

\subsection{Short and long ladders\label{sec-longshort}}

The asymptotic limits of \textit{short ladders} for which $\sqrt \varepsilon L \ll1$ and \textit{long ladders} for which $\sqrt\varepsilon L\gg 1$ are qualitatively different and can be studied in more details. We note that both the ladders shown in Figure~\ref{fig-ladder} and the tri-stable ladder discussed in Section~\ref{sec-stabper} are short since $\sqrt \varepsilon L\simeq  10^{-2} \pi$. 
 
\begin{figure}[h]
\centering
\includegraphics[width=12cm]{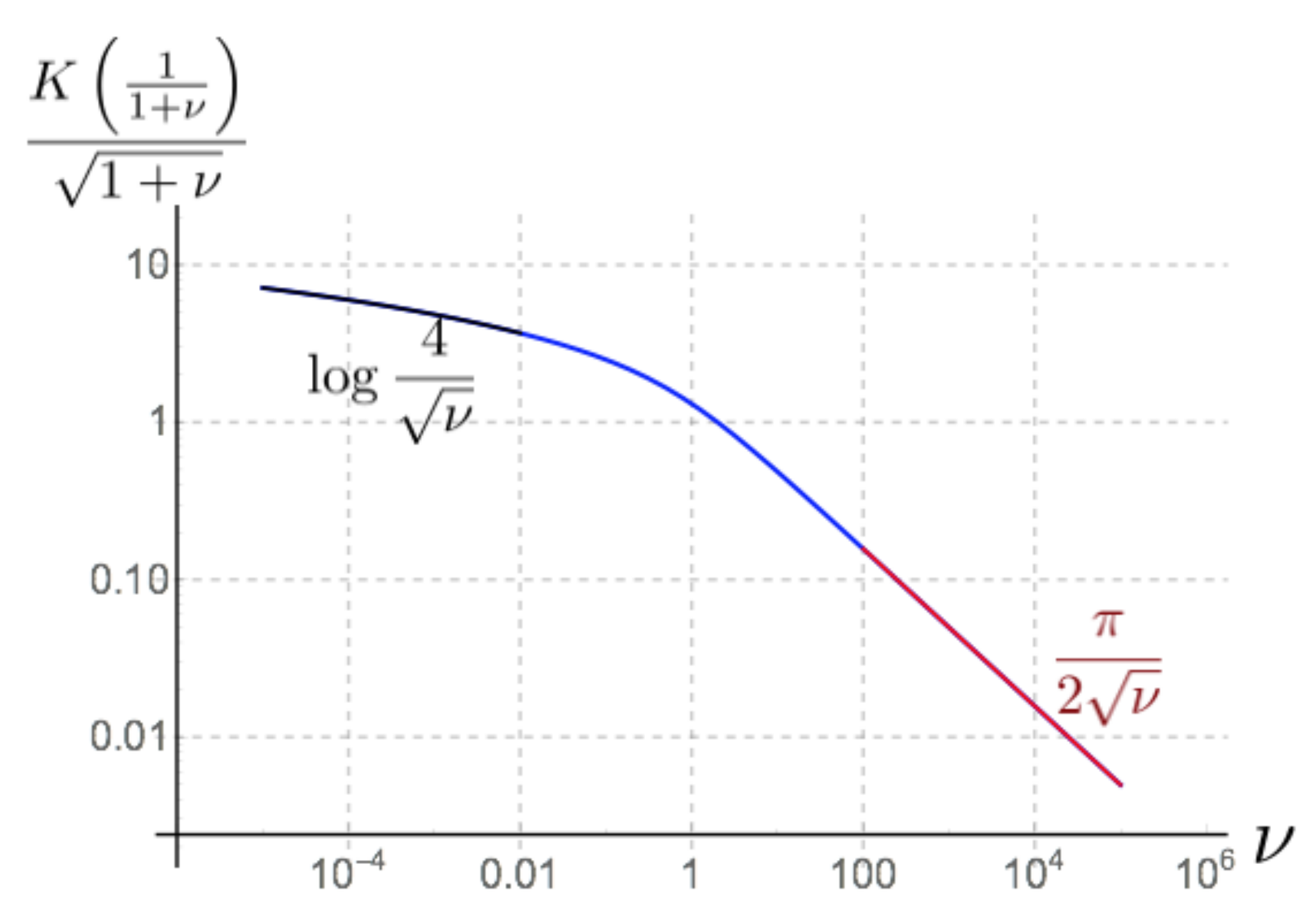}
\caption{A log-log plot of the function $k(\nu)$ defined in\re{ineqper} with its two asymptotic limits for short and long ladders. Note that $k$ is monotonically decreasing. 
}\label{fig-k}
\end{figure}

 The inequalities\re{ineqper} imply that short ladders have perversions with $\nu\gg 1$ and long ladders have perversion for $\nu\ll1$ (see also Figure~\ref{fig-k}). 

\subsubsection*{Short ladders}
For short ladders, we compute 
\be{approxkshort}
k(\nu)= \frac{\pi}{2\sqrt{\nu}} + O\left(\frac 1\nu \right).
\ee
In the limit $\nu\gg1$, the solution of Equation\re{Eqper} is found to be
\be{approxshorteq}
\nu=\hctu_0^2/(4\varepsilon).
\ee
Substituting\re{approxkshort} in\re{ineqper} and\re{approxshorteq} therein gives the bounds on $\hctu_0$ and $L$ so that a short ladder admits a stable perversion for:
\be{boundU0short}
\frac{\pi} L < |\hctu_0| < \frac{2 \pi}{L}. 
\ee

\subsubsection*{Long ladders}
In the limit $\nu\ll1$, the solution of Equation\re{Eqper} is 
\be{sollongeqper}
k(\nu) = \frac 1 2~ \textrm{Arccosh} \left(\frac{2 \varepsilon} { |\hctu_0|} \right)+\frac{\sqrt{\varepsilon} L}{2}. 
\ee
Upon substituting\re{sollongeqper} in\re{ineqper} we find
\be{boundU0long}
\frac{ 2 \sqrt{\varepsilon}} { \cosh(\sqrt{\varepsilon} L/2)} <|\hctu_0| < 2 \sqrt{\varepsilon}. 
\ee

\begin{figure}[h]
\centering
\includegraphics[width=\textwidth]{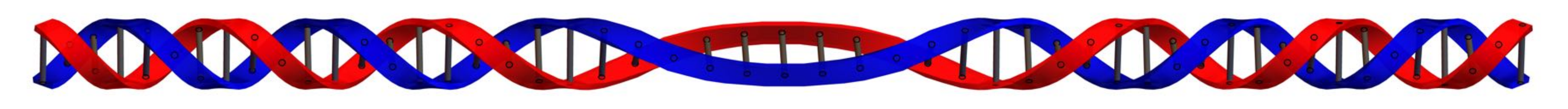}
\caption{Example of a stable perversion of a long ladder: $\varepsilon= 4\, 10^{-2}$, $L=30$, $\hctu_0=0.05$, and $\nu \simeq 0.0024$.  } \label{fig-longper}
\end{figure}

An example of a long-ladder perversion is shown in Figure~\ref{fig-longper}. In this case, stability is guaranteed if $\hctu_0>{ 2 \sqrt{\varepsilon}} /{ \cosh(\sqrt{\varepsilon} L/2)}\simeq 0.040$.

The lower bound in\re{boundU0long} scales like exp$({-\sqrt{\varepsilon}L})$. Accordingly, for long ladders, exponentially small deviation of $\hctu_0$ away from 0 will stabilize the perversion.   

\subsection{Stability of  helical solutions\label{sec-Exacthel}}

Next, we consider the stability of solutions with a given handedness. These equilibria correspond to solutions for $\theta$ centered on either $\pi/2$ or $-\pi/2$ as shown in the bottom-left and bottom-right panels of Figure~\ref{fig-ladder1}. Here, without loss of generality, we consider the left-handed case for which $\theta>0$. This case corresponds to  solutions of the differential equation appearing in\re{ELlad} such that both $\theta(R-L) - \pi/2= -\big( \theta\left( R \right ) - \pi/2\big)$ and $\forall R\in(0,L):\, \theta'(R)\neq0$. 

By repeating the steps of Section~\ref{sec-quantper}, we find 
\be{solthetaLhel}
 \sqrt \varepsilon L~\sqrt{1+\nu} \left (\frac {2 R}{L}-1\right )  = F\left (\theta(R)-\frac\pi 2 \left | -\frac{1}{\nu}\right. \right ). 
\ee

Solutions of\re{solthetaLhel} are stationary functions of\re{LadEnnd} provided that they meet the boundary conditions in\re{ELlad}: 
\be{heleq}
\nu + \textrm{sn}^2\left ( \sqrt{\varepsilon}L \sqrt{\nu} \left | \frac{-1}{\nu}\right. \right) = \frac{\hctu_0^2}{4 \varepsilon}.
\ee
where $ \textrm{sn}(u|m)$ is the Jacobi elliptic sine of parameter $m$.

As long as $\theta(L)\in\left(\frac \pi 2 , \pi \right)$, the left handed helix is stable since  $J[\gamma]=- 1$.
However, if $\theta(L) \in \left(\pi ,\frac{3\pi}{2}\right)$, the solution also intercepts the minima of $V$ at $\theta=0$ and $\theta=\pi$ which implies, $J[\gamma]=1$ and the solution is unstable.  We  use this result to compute the critical value of~\,$\hctu_0$ such that the left-handed helix remains stable. 
The critical case is $\theta(L)=\pi$ which gives the critical value of $\nu$ through\re{solthetaLhel}:
 \be{nucritstaLhel}
\sqrt{\varepsilon} L \sqrt{\nu_\text{crit}} =K(-1/\nu_\text{crit}).
\ee
This  value $\nu_\text{crit}$ 
can be substituted in\re{heleq} to obtain the critical value of $\hctu_0$ such that the helical shape is stable: 
\be{critu0}
\hctu_0= 2 \sqrt{\varepsilon} \sqrt{\nu_\text{crit} }.
\ee
Noting that $\sqrt{\varepsilon} L =K(-1/\nu)/\sqrt{\nu}= k(\nu)$, we can use the asymptotic results of Section~\ref{sec-longshort}.
 For short ladders, we  use the relation $k(\nu)\approx \frac\pi{2 \sqrt{\nu}}$ to find that helical solutions -- centred on $\pm \pi/2$ -- are unstable whenever
\be{limitshort}
|\hctu_0| >  \frac \pi L.
\ee
Similarly, for a {long} ladder ($\sqrt{\varepsilon} L\gg1$), we use $k(\nu)\approx \log \frac 4 {\sqrt{\nu}}$ and conclude that the helical solutions are unstable whenever
\be{limitlong}
|\hctu_0| >  8 \sqrt{\varepsilon} ~~e^{-\sqrt{\varepsilon}L}. 
\ee

\subsection{Effect of $\widehat u$ \label{sec-effectofu}}

So far we have considered the particular case of the functional\re{LadEn} such that   $\widehat u=1$.  We now consider the effect of  $\widehat u\not=0$. In this case the potential  is
%
\be{}
V(\theta)= 4\,b \, (1-\widehat u) \, \cos \theta - (b-\Gamma) \cos 2 \theta.
\ee 
If the flanges are straight  $\hctu_0=\hctu_1=\widehat u=0$,  then $V$ has a different qualitative shape: it presents a single maximum (at $\theta=0$) and no minima on $(0,\pi)$. Accordingly, the only equilibrium of the structure is $\theta=0$ and it is stable. The ladder remains straight.

There is a critical value ${\widehat u}^\star=\frac \Gamma b$ for which $V$ changes qualitatively from having a single maximum (when the straight ladder is the only stable state) to having two maxima (when helical equilibria are stable and there are, typically, multiple equilibria). This transition occurs when $\left .\sdd V \theta \right |_{\theta =0} = 0$: the extremum at $\theta=0$ goes from being a maximum when $\widehat u<\widehat u^\star$ to being a minimum when $\widehat u>\widehat u^\star$. 

The quantitative analysis of Sections~\ref{sec-quantper}-\ref{sec-Exacthel} does not apply when $\widehat u\neq 1$ since, in general,
the solution cannot be expressed in terms of elliptic functions. In this case, the  bounds on $\hctu_0$ for the existence of stable perversions must be found numerically by applying the same criteria.

\section{Conclusions}

Mechanical structures with multiple equilibria represent a rich source of new physical devices. Here, based on the Bristol ladder~\cite{laweda12}, we considered a general class of deployable filamentary structures obtained by connecting two Kirchhoff rods, and showed that they admit many equilibrium solutions. Once a family of equilibria has been identified, the main question is to assess their stability. Here, we used the methods  developed in~\cite{lego15} to probe the stability of these equilibria and to identify regions in parameter space where the system admits three stable states. One of this  state is a perversion consisting of a long ladder with inverted handedness. In particular, we showed how to use the notion of stability index to systematically design stable perversions.

\newpage
\appendix

\section{Birod geometry}\label{app-birod}

\begin{figure}[h]
\centering
\includegraphics[width=6cm]{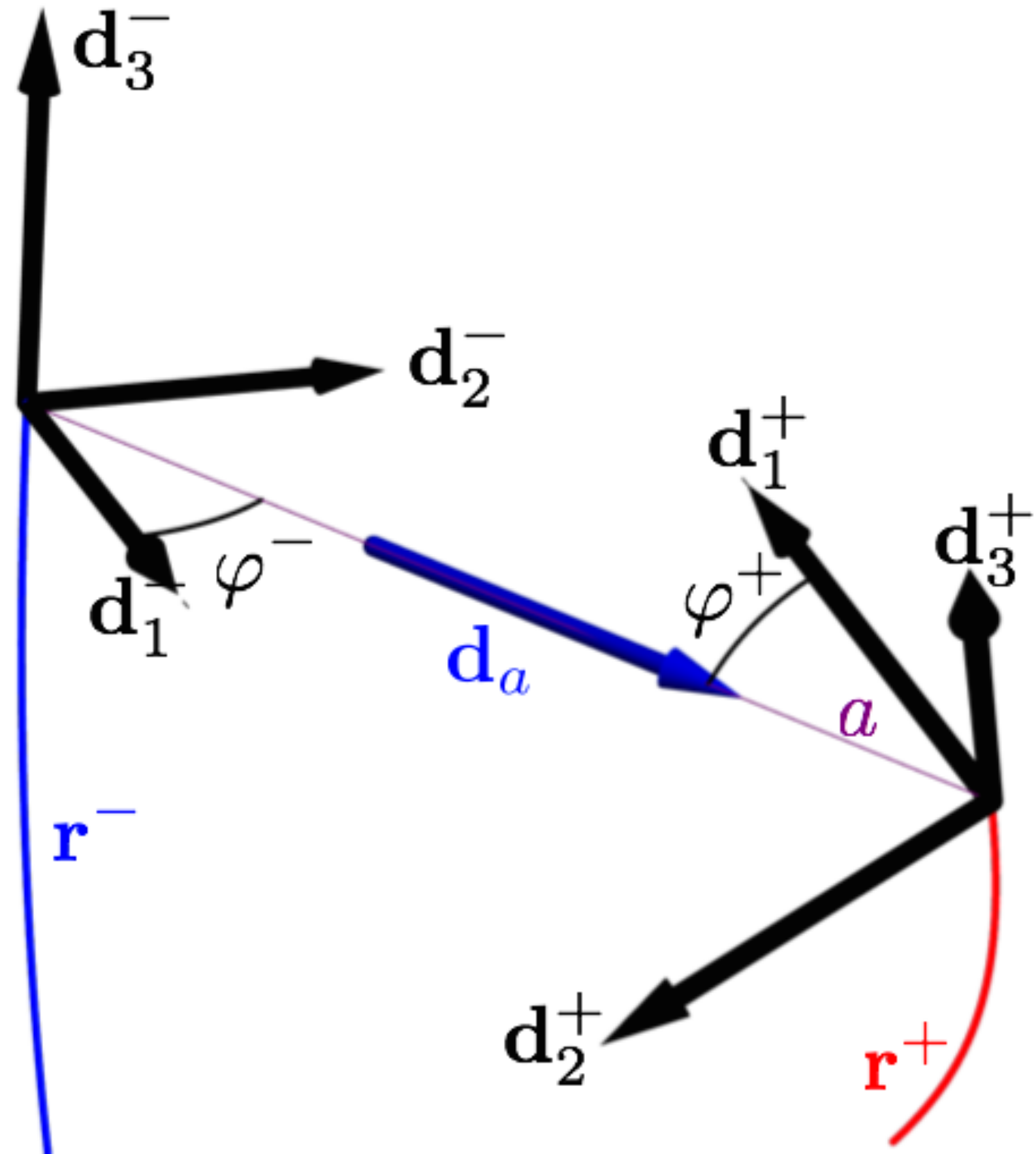}
\caption{Local frames of the sub-rods at centre points $\v r^-$ and $\v r^+$ paired by $G$  together with the angles $\varphi^\pm$, the unit vector along the common cord $\v d_a$ and the distance $a$ between $\v r^-$ and $\v r^+$. } \label{fig-simplerod}
\end{figure}

Following~\cite{lemo15} (see also~\cite{Moakher:2005gcb,Thompson:2002ffa,stva14,Neukirch:2002fka,chmamo09,olbo12}), we define a birod as a pair of rods, the \textit{sub-rods}, attached longitudinally by some mechanism that prevents them from sliding along one another. The sub-rods are labelled by a sign `$-$' or `$+$' and all related variables acquire the corresponding superscript. It is also understood that the sub-rods are parameterised by their respective reference arc-length $S^\pm$.  We furthermore assume that there exists a 1:1 mapping $G:M^-\to M^+$, a number $a$ and two angle functions $\varphi^-(S^-)$ and $\varphi^+(S^+)$ such that in all configurations of the birod (see Figure~\ref{fig-simplerod}), 
\be{assume}
\left \{
\begin{split}
\v r^+ \left(G(S^-)\right)&= \v r^-\left( S^-\right) + a\left (\cos \varphi^-\, \v d_1^- + \sin \varphi^-\, \v d_2^- \right) (S^-), \\
\v r^- \left(G^{-1}(S^+)\right)&= \v r^+\left( S^+\right) + a \left(\cos \varphi^+\,\v d_1^+ + \sin \varphi^+\, \v d_2^+ \right) (S^+).
\end{split}
\right .
\ee
The assumption\re{assume} is quite restrictive: it implies that the distance $a$ between centre points of sections connected by $G$ is constant along the length of the structure and unchanging in all configurations -- a rigid constraint. Furthermore, the point $\v r^+(G)$ (resp. $\v r^-(G^{-1})$) must be in the plane containing $\v r^-$ (resp. $\v r^+$) and spanned by $\v d_1^-$ and $\v d_2^-$ (resp. $\v d_1^+$ and $\v d_2^+$). It will be useful to have a definition of the unit vector $\v d_a$ along the direction of the common cord $\v r^+(G)-\v r^-$:
\be{defda}
\v d_a= \cos \varphi^- \, \v d_1^- + \sin\varphi^- \, \v d_2^-=-( \cos \varphi^+ \, \v d_1^+ + \sin\varphi^+ \, \v d_2^+).
\ee
Note in particular that $\v d_a$ is perpendicular to both $\v d_3^-$ and $\v d_3^+$. We define the angle $\theta$ between the tangents to the sub rods $\v d_3^\pm$ such that a rotation of angle $\theta$ about the axis $\v d_a$ brings $\v d_3^-$ on $\v d_3^+$. 

With these definitions, a birod can be viewed as a one dimensional manifold $M$ defined as a collection of pairs of sections. Each point of $M$ contains one point of $M^-$ and one points of $M^+$ the manifolds describing the sub-rods. Let $R$ be a material parameterisation of $M$, we define the 1:1 mappings $G^\pm:R\to S^\pm$. The derivatives of these mappings will prove useful so we define 
\be{defg}
g^-=\sd{G^-}{R},\qquad g^+ = \sd{G^+}{R},\quad \textrm{and}\quad g=\sd{G}{S^-} \qquad \textrm{such that } \quad g=g^+/g^-.
\ee

It is then possible to prove (see~\cite{lemo15}) that Equations\re{assume} imply
\be{defU}
 g^-\v u^- + (\varphi^-)' \, \v d_3^- + \frac {\theta'} 2 \v d_a= g^+\v u^+ + (\varphi^+)' \, \v d_3^+ - \frac {\theta'} 2 \v d_a =: \vtu,
\ee
where $'$ denotes derivation w.r.t. $R$.

Furthermore, there exists a function $\ctu(R)$ such that (see~\cite{lemo15} for the proof)
\be{Ual}
\vtu =\ctu\,  \v d_a + \frac 1 a\, \v d_a \times \left (\alpha^+ g^+ \v d_3^+ - \alpha^- g^- \v d_3^-\right),
\ee
where $\alpha^\pm=\sd{s^\pm}{S^\pm}$ are the stretches of the sub-rods. 

The vector $\vtu$ is to a birod what $\v u$ is to a Kirchhoff rod. It expresses the rate of rotation of the non-orthogonal frame $\{\v d_a,\v d_3^-,\v d_3^+\}$ per unit $R$:
\be{rotframe}
{\v d_3^\pm}'=\left  (\vtu \pm \frac{\theta'}{2} \v d_a \right ) \times \v d_3^\pm, \quad \textrm { and }\quad {\v d_a}'=\vtu \times \v d_a.
\ee

Accordingly, given the four functions $\theta'$, $\alpha^+$, $\alpha^-$ and $\ctu$ and initial values for $\v d_3^+$, $\v d_3^-$ and $\v d_a$ at some particular value of $R$, we can substitute $\vtu$ according to\re{Ual} in the system of 9 first-order differential equations\re{rotframe} to obtain an initial value problem. A direct (numerical) integration reveals the shape of the structure. 

Alternativelay, we can substitute\re{Ual} in\re{defU} to obtain a direct expression for the curvatures of the sub-rods:
 \be{upmUalth}
\v u^\pm = \frac {1}{ a} \left ( \frac {a}{g^\pm}\left [ \ctu  \pm \frac{\theta'}{2} \right ] \v d_a+ \frac {g^-}{ g^\pm} \v d_a \times \left ( \alpha^+ g\, \v d_3^+ -\alpha^- \v d_3^-\right) - \frac{(\varphi^\pm)'}{  g^\pm}  \v d_3^\pm\right ).
\ee
Then it is possible to integrate\re{shape1rod} for each sub-rod to find their respective shape.

\section{Potential energy of the birod\label{app-energy}}

We compute energy stored in a birod due to the bending of its sub-rod and express it as a functional of $\alpha^\pm$, $\ctu$, $\theta$ and $\theta'$. We assume that the sub-rods are hyperelastic and have quadratic bending-energy densities:
\be{}
W_{\textrm{bend}}^\pm= 
\frac{(EI)^\pm}{2}\left [ (\c u_1 - \widehat{\c u}_1)^2+ 
b^\pm (\c u_2 - \widehat{\c u}_2)^2+
\Gamma^\pm (\c u_3 - \widehat{\c u}_3)^2 \right].
\ee
where $(EI)^\pm$, $b^\pm$ and $\Gamma^\pm$ are material parameters.

The bending energy stored in the birod is then defined as 
\ben{appDefEbend}
\mathscr E_{\textrm{bend}}&=& \int_a^{b} \textrm d R \Big\{ g^+ W_{\textrm{bend}}^+ + g^- W_{\textrm{bend}}^- \Big \},\\
&=&\int_a^{b} \textrm d R\frac{ (EI) g^-}{2} \Bigg\{  g^2 \mathcal E^+ \bigg( (\c u_1^+ - \widehat{\c u}_1^+)^2 + b^+ (\c u_2^+ - \widehat{\c u}_2^+)^2+ \Gamma^+(\c u_3^+ - \widehat{\c u}_3^+)^2 \bigg)  \nonumber\\
&&\qquad\qquad\qquad\qquad + \mathcal E^- 
 \bigg( (\c u_1^- - \widehat{\c u}_1^-)^2 + b^- (\c u_2^- - \widehat{\c u}_2^-)^2  + \Gamma^-(\c u_3^- - \widehat{\c u}_3^-)^2 \bigg)\Bigg \},\label{appEbend}
 \een
 where $(EI)=\frac{(EI)^+}{g} + (EI)^-$, $\mathcal E^+ = \frac{(EI)^+}{g (EI)}$ and $\mathcal E^- =  \frac{(EI)^-}{(EI)}$ so that $\mathcal E^+ + \mathcal E^-=1$.
 
Next, recall\re{upmUalth} from Appendix~\ref{app-birod}:
 \be{appupmUalth}
\v u^\pm = \frac {1}{ a} \left ( \frac {a}{g^\pm}\left [ \ctu  \pm \frac{\theta'}{2} \right ] \v d_a+ \frac {g^-}{ g^\pm} \v d_a \times \left ( \alpha^+ g\, \v d_3^+ -\alpha^- \v d_3^-\right) - \frac{(\varphi^\pm)'}{  g^\pm}  \v d_3^\pm\right ).
\ee
 Also, since $\v d_a$ is in both the planes $(\v d_1^-,\v d_2^-)$ and $(\v d_1^+,\v d_2^+)$  respectively perpendicular to $\v d_3^-$ and $\v d_3^+$, we have 
 \ben{rotdmda}
 &&\left \{
 \begin{split}
 \v d_1^-&= \cos\varphi^- \v d_a - \sin \varphi^- \v d_3^-\times \v d_a,\\
 \v d_2^-&=\sin\varphi^- \v d_a + \cos \varphi^- \v d_3^-\times \v d_a,
 \end{split} 
 \right .\\
 \textrm{and}&&
 \left \{
 \begin{split}
  \v d_1^+&= \cos(\varphi^+ + \pi) \v d_a - \sin (\varphi^+ + \pi) \v d_3^+\times \v d_a,\\
 \v d_2^+&=\sin(\varphi^+ + \pi) \v d_a + \cos (\varphi^+ + \pi) \v d_3^+\times \v d_a.
 \end{split} 
 \right .\label{rotdpda}
 \een

Substituting~(\ref{appupmUalth},\ref{rotdmda},\ref{rotdpda}) in $\c u_{1,2,3}^\pm= \v u^\pm\cdot \v d_{1,2,3}^\pm$, we find 
\ben{u1pmalpmtheta}
\c u_1^\pm &=&- \frac {1}{g^\pm} \cos \varphi^\pm \left [  \frac {\theta'}{2}\pm \ctu   \right ] -\frac {1}{ g^\pm} \sin \varphi^\pm \frac{ \alpha^\pm g^\pm - \cos \theta \alpha^\mp  g^\mp}{a},\\
\c u_2^\pm &=&- \frac {1}{ g^\pm} \sin \varphi^\pm \left [ \frac {\theta'}{2}\pm \ctu\right ] +\frac {1}{ g^\pm} \cos \varphi^\pm \frac{\alpha^\pm g^\pm - \cos \theta \alpha^\mp  g^\mp}{a},\label{u2pmalpmtheta}\\
\c u_3^\pm &=& - \frac {1 }{g^\pm} \frac{\alpha^\mp g^\mp}{a} \sin \theta -\frac{(\varphi^\pm)'}{g^\pm}.\label{u3pmalpmtheta}
\een

Hence, by direct application of\re{u3pmalpmtheta},
\be{u3mu3hatsq}
\begin{split}
(\c u_3^+ -\widehat{\c u}_3^+)^2&= \frac 1 {{g^+}^2} \left( \frac{g^-}{a} \,\alpha^- \sin\theta +  \left ( g^+ \widehat{\c u}_3^+ +(\varphi^+)' \right )\right )^2\\
(\c u_3^- -\widehat{\c u}_3^-)^2&=\frac1 {{g^- }^2} \left(\frac{g^-}{a}\, \alpha^+ g \sin\theta +  \left (g^- \widehat{\c u}_3^- + \varphi^-)' \right )\right )^2
\end{split}
\ee 

Next we define
\be{Defhctu}
\begin{split}
&\hctu_1^\pm  =\mp g^\pm (\widehat{\v u}^\pm \cdot \v d_a)= g^\pm~ \left( \widehat{\c u}_1^\pm \cos\varphi^\pm +\widehat{\c u}_2^\pm \sin\varphi^\pm\right ),\\
&\hctu_2^\pm = \mp  g^\pm \, (\widehat{\v u}^\pm \cdot (\v d_3^\pm\times \v d_a))= g^\pm~ \left( -\widehat{\c u}_1^\pm \sin\varphi^\pm +\widehat{\c u}_2^\pm \cos\varphi^\pm\right ) ,
\end{split}
\ee
which can be inverted to obtain 
\be{Defhctuinvert}
\begin{split}
\widehat{\c u}_1^\pm &= \frac {1}{g^\pm} \left( \cos \varphi^\pm \hctu_1^\pm - \sin \varphi^\pm \hctu_2^\pm\right ), \\
\widehat{\c u}_2^\pm &= \frac {1}{g^\pm} \left( \sin \varphi^\pm \hctu_1^\pm + \cos \varphi^\pm \hctu_2^\pm\right ). 
\end{split}
\ee
Gathering~(\ref{u1pmalpmtheta},\ref{u2pmalpmtheta},\ref{Defhctuinvert}) yields
\be{uimuihatp}
\begin{split}
\c u_1^+- \widehat{\c u}_1^+ &=\frac{-1}{g^+} \bigg [ \cos \varphi^+ \Big(\theta'/2 + \ctu +\hctu_1^+ \Big) + \sin \varphi^+ \Big(  \frac{g^-}{a} (\alpha^+ g - \alpha^- \cos \theta) -\hctu_2^+ \Big) \bigg],\\
\c u_2^+- \widehat{\c u}_2^+ &=\frac{1}{g^+} \bigg [ -\sin \varphi^+ \Big(\theta'/2 + \ctu +\hctu_1^+\Big) + \cos \varphi^+ \Big( \frac{g^-}{a}(\alpha^+ g - \alpha^- \cos \theta) -\hctu_2^+\Big) \bigg].
\end{split}
\ee

Next define $\rho^+$ and $\omega^+$ such that
\be{Defrhoomegap}
\rho^+ \cos \omega^+ =\frac{\theta'}{2} + \ctu +\hctu_1^+, \qquad \textrm{and} \qquad 
\rho^+ \sin \omega^+ =(\alpha^+ g - \alpha^- \cos \theta) \frac{g^-}{a} -\hctu_2^+.
\ee

Using the definitions\re{Defrhoomegap} and substituting according to\re{uimuihatp},
\ben{}
&&\big(\c u_1^+- \widehat{\c u}_1^+\big)^2+ b^+ \big(\c u_2^+- \widehat{\c u}_2^+\big)^2\nonumber \\
&&\qquad =\left (\frac{\rho^+}{g^+}\right)^2 \bigg[ \big(\cos\varphi^+ \cos \omega^+ +\sin\varphi^+ \sin\omega^+\big)^2 + b^+\big (-\sin\varphi^+ \cos \omega^+ + \cos\varphi^+ \sin\omega^+\big)^2\bigg]\nonumber\\
&&\qquad =\left (\frac{\rho^+}{g^+}\right)^2 \bigg[ \frac{1+b^+}{2} +\frac{1-b^+}{2} \cos^2(2 \omega^+ - \varphi^+)]\nonumber\\
&&\qquad =\left(\frac{1}{g^+}\right)^2 \bigg[ \frac{1+b^+}{2} {\rho^+}^2   + \frac{1-b^+}{2} \Big(\cos2\varphi^+  {\rho^+}^2(\cos^2\omega^+ -\sin^2 \omega^-) + 2 \sin 2 \varphi^+ {\rho^+}^2 \cos\omega^+ \sin\omega^+ \Big) \bigg]\nonumber\\
&&\qquad = \left(\frac{1}{g^+}\right )^2 \left[\frac{1+b^+}{2}\left (\left(\frac{\theta'}{2} + \ctu +\hctu_1^+ \right)^2+\left(\frac{g^-}{a} (\alpha^+ g - \alpha^- \cos \theta) -\hctu_2^+\right)^2 \right)  \right . \nonumber \\
&& \qquad\qquad \qquad + \frac{1-b^+}{2} \Bigg( \cos 2\varphi^+ \left( \left(\frac{\theta'}{2} + \ctu +\hctu_1^+\right)^2 -\left(\frac{g^-}{a}(\alpha^+ g - \alpha^- \cos \theta) -\hctu_2^+\right)^2 \right) \nonumber\\
&&\left.\qquad\qquad \qquad \qquad \qquad +2 \sin2 \varphi^+\left(\frac{\theta'}{2} + \ctu +\hctu_1^+\right) \left(\frac{g^-}{a}(\alpha^+ g - \alpha^- \cos \theta )-\hctu_2^+\Bigg)
\right )
\right ] .\label{Wp12}
\een

A similar argument leads to 
\ben{}
&&\big(\c u_1^- - \widehat{\c u}_1^-\big)^2+ b^- \big(\c u_2^- - \widehat{\c u}_2^-\big)^2\nonumber \\
&&\qquad = \left(\frac{1}{g^-}\right)^2 \left[\frac{1+b^-}{2}\left (\left(\frac{\theta'}{2} - \ctu +\hctu_1^-\right)^2+\left(\frac{g^-}{a}(\alpha^- - \alpha^+ g \cos \theta )-\hctu_2^-\right)^2 \right)  \right . \nonumber \\
&& \qquad\qquad \qquad + \frac{1-b^-}{2} \Bigg( \cos 2\varphi^- \left( \left(\frac{\theta'}{2} - \ctu +\hctu_1^-\right)^2 -\left(\frac{g^-}{a}(\alpha^- - \alpha^+ g \cos \theta)-\hctu_2^-\right)^2 \right)  \nonumber\\
&&\left.\qquad\qquad  \qquad\qquad \qquad +2 \sin2 \varphi^-\left(\frac{\theta'}{2} - \ctu +\hctu_1^-\right) \left(\frac{g^-}{a}(\alpha^- - \alpha^+ g \cos \theta) -\hctu_2^-\Bigg)
\right )
\right ] .\label{Wm12}
\een

While\ret{Wp12}{Wm12} are valid for sub-rods with general cross-sections, we focus here on two special cases. Indeed, since we have assumed that the distance $a$ between the centrelines of the sub rods is constant, if the cross-sections of the sub rods are not symmetric ($b^\pm\neq1$), we assume first that $\varphi^\pm = k^\pm \pi$ (with $k^\pm \in \mathbb Z$). If the sub-rod $\chi$ (with $\chi=+$ and/or $-$) is symmetric: $b^\chi=1$, then we allow for general function(s) $\varphi^\chi(S^\chi)$. The assumption does lead to a (slight) loss of generality but birods with constant distance $a$ which does not respect it are seldom in practice. Let us mention two such examples however: an asymmetric rod encompassed in a slender matrix and a birod the sub-rod of which have circular cross-sections with anisotropic material properties. For such rods, all three terms must be kept in\ret{Wp12}{Wm12} and due to the mixed term (multiplied by $\sin 2 \varphi$) there is a coupling between the rate of curvature $\ctu$ around $\v d_a$ and the deformations. Here for the sake of simplicity, we assume this is not the case and either
\begin{itemize}
\item{the cross-sections are symmetric so that $b^\pm=1$ or,}
\item{the common chord is aligned with the principal direction $\v d_1^\chi$ in which case $\cos 2 \varphi^\chi=1$ and $\sin 2 \varphi^\chi=0$.}
\end{itemize}

In both those cases,\ret{Wp12}{Wm12} become 
\be{Wp12simp}
\begin{split}
\big(\c u_1^+- \widehat{\c u}_1^+\big)^2+ b^+ \big(\c u_2^+- \widehat{\c u}_2^+\big)^2&= \left(\frac{1}{ g^+}\right)^2 \left[\left(\frac{\theta'-\hctu_0}{2} + \left(\ctu -\frac {\hctu_1}{2} \right ) \right)^2+ b^+ \left(\frac{g^-}{a}(\alpha^+ g - \alpha^- \cos \theta)-\hctu_2^+\right)^2 \right ] ,\\
\big(\c u_1^- - \widehat{\c u}_1^-\big)^2+ b^- \big(\c u_2^- - \widehat{\c u}_2^-\big)^2&= \left(\frac{1}{g^-}\right)^2 \left[\left(\frac{\theta'-\hctu_0}{2} - \left(\ctu -\frac {\hctu_1}{2} \right ) \right)^2+ b^- \left(\frac{g^-}{a}(\alpha^- - \alpha^+ g \cos \theta) -\hctu_2^-\right)^2 \right],
\end{split}
\ee
where we defined $\hctu_0=-\hctu_1^+ - \hctu_1^-$ and $\hctu_1=-\hctu_1^+ +\hctu_1^-$. Note that using\re{Defhctu} in these definitions, we obtain $\hctu_0=(g^+ \widehat{\v u}^+- g^-\widehat{\v u}^-) \cdot \v d_a$ and $\hctu_1= (g^+ \widehat{\v u}^++ g^-\widehat{\v u}^-) \cdot \v d_a$. 

Substituting~(\ref{u3mu3hatsq},\ref{Wp12simp}) in\re{appEbend}, we obtain 
\ben{}
\mathscr E_{\textrm{bend}}& =&\int_a^{b}\textrm{d} R\frac {(EI)}{2 g^-}  \Bigg\{
 \left (\frac{\theta'-\hctu_0}{2}\right )^2 + \left (\ctu-\frac{\hctu_1}{2}\right )^2 + 2 (\mathcal E^+- \mathcal E^-) \left (\frac{\theta'-\hctu_0}{2} \right ) \left ( \ctu -\frac {\hctu_1}{2} \right )\nonumber  \\
&&\qquad +b^-\, \mathcal E^- \left [\frac{g^-}{a}(\alpha^- -\alpha^+ g\cos \theta ) - \hctu_2^- \right ]^2 
+b^+\, \mathcal E^+ \left [\frac{g^-}{a}(\alpha^+ g - \alpha^- \cos \theta) - \hctu_2^+ \right ]^2 \label{appEbendfin}\\
&& \qquad + \Gamma^-\mathcal E^- \left [ \frac{g^-}{a}\, \alpha^+ g  \sin\theta  +  \left (g^- \widehat{\c u}_3^- + (\varphi^-)'\right )\right ]^2
+ \Gamma^+ \mathcal E^+ \left [\frac{g^-}{a}\, \alpha^- \sin \theta + \left (g^+ \widehat{\c u}_3^+ + (\varphi^+)' \right ) \right ]^2 \Bigg\} .\nonumber
\een

If we furthermore assume that the sub-rods are inextensible, Equation\re{appEbendfin} can be compactly expressed as 
\ben{appRefEnergyInext}
&&\qquad\qquad\qquad \mathscr E=  \int_a^{b}\frac{ (EI)}{2 g^-} ~ \bigg\{\mathcal L[\ctu,\theta,R] +C(R) \bigg\} ~\textrm{d}R\\
&&\textrm{with} \quad \mathcal L= 
\left (\frac{\theta'-\hctu_0}{2}\right )^2 +  \left ( \, \ctu-\frac{\hctu_1}{2 }\right )^2+ 2\left(\mathcal E^+ - \mathcal E^-\right)\frac{\theta'-\hctu_0 }{2} \left ( \ctu -\frac {\hctu_1}{2} \right )- 2V(\theta),\nonumber\\
&&\textrm{and}\qquad V(\theta)= \cos \theta ~\left (  \left(\frac{g^-}{a}\right)^2 B - \frac{g^-}{a} \hctu_2 \right) - \frac{g^-}{a}\, \frac{A}{4}\,  \cos 2 \theta -\frac{g^-}{a}\, \sin\theta\, \hctu_3 - \frac{\psi(\theta)}{2},\nonumber
\een
where $\hctu_2$, $\hctu_3$, $A$, $B$ and $C$ are functions of $R$ defined by 
\ben{}
\hctu_{2}&=&b^+ \mathcal E^+ \hctu_{2}^+ + g\, b^- \mathcal E^-\, \hctu_{2}^- ,\nonumber\\
\hctu_{3}&=&\Gamma^+  \mathcal E^+ \hctu_{3}^+ +g\, \Gamma^- \mathcal E^- \hctu_{3}^- \nonumber\\
A&=&g^2 \mathcal E^-  (b^--\Gamma^-) + \mathcal E^+ (b^+ - \Gamma^+),\label{defconst}\\
B&=&g (b^- \mathcal E^-+b^+ \mathcal E^+), \nonumber\\
C&=& b^- \mathcal E^-\left ( \frac{g^-}{a}-\hctu_2^-\right )^2 + b^+ \mathcal E^+ \left ( \frac{g^-}{a}\, g - \hctu_2^+\right)^2+\Gamma^- \mathcal E^- \, (\hctu_3^-)^2+\Gamma^+ \mathcal E^+\, (\hctu_3^+)^2\nonumber \\
&&\qquad\qquad +\frac{g^2 (b^- + \Gamma^-) \mathcal E^-+ (b^+ +\Gamma^+) \mathcal E^+}{2},\nonumber
\een
where $\hctu_{3}^\pm=g^\pm \widehat {\c u}_3^\pm + (\varphi^\pm)'$. Note that the last term in\re{appRefEnergyInext} may be safely omitted since the potential energy of a structure is defined up to an arbitrary additive constant. 

Finally, with the parameters chosen in the main text, we find $(EI)^+=(EI)^-=(EI)\, a$, $g^\pm=a$, $g=1$, $b^+ =b^-=b$, $\mathcal E^+ = \mathcal E^- = 1/2$, $\Gamma^+ = \Gamma^-= \Gamma$, $A=b-\Gamma$, $B=b$, $\hctu_2= b \, \widehat u$, and $\hctu_3=0$ which can be substituted in\re{appRefEnergyInext} to give\re{LadEn}.

 \section{Comparing different models}\label{app-originalmodel}

If the flanges have narrow rectangular cross sections of width $w$ and thickness $h$ with $h\ll w$, their bending stiffness can be estimated as $(EI)_i = E\,  I_i$ where $E$ is the Young modulus of the material and $I_i$ the second moment of area of the section about the director $\v d_i$: $I_1= w^2 \frac{w h} {12}$ and $I_2= h^2 \frac{wh}{12}$. In that case, the parameter $b=(h/w)^2\ll1$. The torsional stiffness $(\mu J)=\Gamma\, (EI)$ of such a bar can be approximated  by $(\mu J) = 4\, \mu \, h^2 \frac{wh}{12}$  where $\mu$ is the shear modulus of the material~\cite{ti50TorStif}. Accordingly, the parameter $\Gamma = 4 b \frac{\mu}{E}$. Substituting these estimates in\re{LadEn}, the energy of the ladder becomes
\be{LadEnEst}
\mathscr E= \frac{ E w^3 h}{48} \int_0^{L}   
\frac{\left (\theta'-\hctu_0\right )^2 }{2}+ 2 \left ( \, \ctu-\frac{\hctu_1}{2 }\right )^2-  4 \left (\frac{h}{w}\right)^2 \bigg( \cos \theta ~\left( 1 - \widehat{u} \right) - \left (\frac{1}{4}-\frac{\mu}{E}\right) \,  \cos 2 \theta \bigg)~ \textrm{d} R.
\ee

In~\cite{laweda12} the flanges are modeled as thin inextensible shells. An important consequence of that modeling choice is that such shells are unable to bend about the direction $\v d_1^\pm$ (see Figure~\ref{fig-ladder}): $\c u_1^\pm=0$. Substituting $\varphi^\pm= \left(\frac 1 2 \pm \frac 1 2\right) \pi$ and $g^\pm=a$ in Equation\re{u1pmalpmtheta} of Appendix~\ref{app-energy}, we find
\be{}
\c u_1^\pm= \frac 1 a  \left (\ctu \pm \frac{\theta'}{2}\right ). 
\ee
Hence the assumption $\c u_1^\pm=0$ leads to the constraints $\theta'=\ctu=0$ in our model. In that case, the energy of the ladder\re{LadEnEst}, simply becomes 
\be{approxE}
\mathscr E_{\textrm{approx}} = \frac{E w h^3 L^-}{12 a} \left [ 
\left(\frac 1 4 - \frac \mu E \right ) \cos  2 \theta - \cos \theta~ (1-  \widehat { u} ) 
\right ],
\ee
where $\theta$ is a constant. For comparison,  the energy $U$ found in~\cite{laweda12} for flanges constituted  of  a $[0_5]$ lay-up is\footnote{See Equation (4.1) in~\cite{laweda12} with $D_{16}^\star=0$ for the $[0_5]$ lay-up. Note that their $\theta$ is half of ours and their parameter $\alpha =R_i/R$ is the ratio between $R_i$ the radius of curvature of the flanges and  half-length $R$  of the spokes: $\alpha= 2/(a \widehat{\c u})$.}: 
\be{defULaweda}
U=2 D_{11}^\star \frac{L^- w} {a^2}  \Big[ 
\left (  \frac{1}{4} - \frac{ D_{66}^\star}{D_{11}^\star} \right ) \cos 2 \theta - \cos \theta~ (1- a\widehat{\c u} ) + \overbrace{\frac{ D_{66}^\star}{D_{11}^\star} + \frac 1 2  (a\widehat{\c u}-1)^2}^{\textrm{a constant independant of }\theta}
\Big ],
\ee
where $D_{11}^\star$ and $D_{66}^\star$ are specific components of the reduced flexural stiffness matrix of the flange as defined from lamination theory and $a \,\widehat{\c u}= \widehat u$: that is $\widehat{\c u}$ is the dimensional curvature corresponding to the non-dimensional $\widehat u$. For the prototype realized in~\cite{laweda12}, the reported values were $D_{11}^\star= 1.986~\textrm{Nm}$ and $D_{66}^\star=0.0603~\textrm{Nm}$ together with the length of the spokes $a=57 ~\textrm{mm}$, the reference curvatures $\widehat{\c u} = 1/57 ~\textrm{mm}^{-1}$ and the flanges length $L^+=L^-=179~\textrm{mm}\simeq \pi a$, width $w=9~\textrm{mm}$ and thickness $h=0.11~\textrm{mm}$.  Since the energy of the system is defined up to an additive constant, the last two terms in\re{defULaweda} are inconsequential. The comparison between\re{approxE} and\re{defULaweda} allows us to estimate the key parameters of our model for the purpose of application to their prototype: $\mu/E=D_{66}^\star/D_{11}^\star\simeq 0.030$ (this is to be understood as a measure of the effective shear modulus of the flanges) and $\sqrt b = h/w = 0.11/9\simeq 0.012$ which gives $\varepsilon = b-\Gamma= b\, (1-4\, \mu/E)\simeq1.27\, 10^{-4}$. We also compute $\sqrt{\varepsilon} L \simeq 0.02\ll 1$ so the prototype of the Bristol ladder is a short ladder.

Perversions were not observed in~\cite{laweda12} for two reasons. First, with the prototype's parameters a perversion would require the flanges to be longer than $L_0\,  a \simeq 7.9 \, \textrm{m}$. For comparison, the prototype had flanges of length $L^-=   0.179\, \textrm{m}$. Second, even if a prototype of sufficient length was built with all other parameters unchanged, all perverted equilibria would be unstable since $\hctu_0=0$. The analysis of Section~\ref{sec-quantit} gives an effective way to build a new prototype with stable perverted states.

 \end{document}